\documentclass[thmsa]{article}
\usepackage{amsfonts}
\usepackage{amssymb}
\usepackage{amsmath}

\input{tcilatex}
\begin{document}

\begin{center}
{\Large Fractals of generalized F}$-${\Large \ Hutchinson operator in
b-metric spaces}

\textbf{Talat Nazir}$^{1,2},$ \textbf{Sergei Silvestrov}$^{1}$ and \textbf{%
Xiaomin Qi}$^{1}$

{\small $^{(1)}$\textit{Division of Applied Mathematics, School of
Education,\ Culture and\ Communication,}} {\small \textit{M\"{a}lardalen
University, 72123 V\"{a}ster\aa s, Sweden}}

$^{(2)}${\small \textit{Department of Mathematics, COMSATS Institute of
Information Technology, 22060 Abbottabad, Pakistan\medskip }}

E-mail{\small : talat@ciit.net.pk, sergei.silvestrov@mdh.se,
xiaomin.qi@mdh.se}\bigskip
\end{center}

\noindent
--------------------------------------------------------------------------------------------

\noindent \textit{Abstract:} \ \ The aim of this paper is to construct a
fractal with the help of a finite family of generalized $F$-contraction
mappings, a class of mappings more general than contraction mappings,
defined in the setup of $b$-metric space. Consequently, we obtain a variety
of results for iterated function system satisfying a different set of
contractive conditions. Our results unify, generalize and extend various
results in the existing literature.

\noindent \textbf{---------------------------------------------}


\noindent \textit{2000 \ Mathematics Subject Classification: }\texttt{47H10,
54E50, 54H25}.

\noindent \textbf{---------------------------------------------}\medskip

\section{Introduction and preliminaries}

Iterated function systems are method of constructing fractals and are based
on the mathematical foundations laid by Hutchinson \cite{Hutchinson81}. He
showed that Hutchinson operator constructed with the help of a finite system
of contraction mappings defined on a Euclidean space $%
\mathbb{R}
^{n}$ has closed and bounded subset of $%
\mathbb{R}
^{n}$ as its fixed point, called attractor of iterated function system (see
also in \cite{Barnsley93} ). In this context, fixed point theory plays
significant and vital role to help in construction of fractals.

Fixed point theory is studied in an environment created with appropriate
mappings satisfying certain conditions. Recently, many researchers have
obtained fixed point results for single and multi-valued mappings defined on
metrics spaces.$\ $Banach contraction principle \cite{Banach} is of
paramount importance in metrical fixed pint theory with a wide range of
applications, including iterative methods for solving linear, nonlinear,
differential, integral, and difference equations. This initiated several
researchers to extend and enhance the scope of metric fixed point theory. As
a result, Banach contraction principle have been extended either by
generalizing the domain of the mapping \cite{Abdeljawad11, AKR11, HZ07,
IAN15, KRR09, Tarafdar74} or by extending the contractive condition on the
mappings \cite{BW69, Edelstein, Kirk03, Meir69, Rakotch}. There are certain
cases when the range $X$ of a mapping is replaced \ with a family of sets
possessing some topological structure and consequently a single valued
mapping is replaced with a multivalued mapping.\ Nadler \cite{Nadler69} was
the first who combined the ideas of multivalued mappings and contractions
and hence initiated the study of metric fixed point theory of multivalued
operators, see also \cite{AAN15, AssadKirk, SVetro13}. The fixed point
theory of multivalued operators provides important tools and techniques to
solve the problems of pure, applied and computational mathematics which can
be re structured as an inclusion equation for an appropriate multivalued
operator.

The concept of metric has been generalized further in one to many ways. The
concept of a $b$-metric space was introduced by Czerwik in \cite{CZ}. Since
then, several papers have been published on the fixed point theory of
various classes of single-valued and multi-valued operators in $b$-metric
space \cite{AN, AY, BO, BO1, C4, CH, CZ, CZ1, CZ2, KU, RO}.

In this paper, we construct a fractal set of iterated function system, a
certain finite collection of mappings defined on a $b$-metric space which
induce compact valued mappings defined on a family of compact subsets of a $%
b $-metric space. We prove that Hutchinson operator defined with the help of
a finite family of generalized $F$-contraction mappings on a complete $b$%
-metric space is itself generalized $F$-contraction mapping on a family of
compact subsets of $X.$ We then obtain a final fractal obtained by
successive application of a generalized $F$-Hutchinson operator in $b$%
-metric space.\medskip

\noindent \textbf{Definition 1.1.} \ \ Let $X$ be a nonempty set and $b\geq
1 $ a given real number. A function $d:X\times X\rightarrow 
\mathbb{R}
^{+}$ is said to be a $b$-metric if for any $x,y,z\in X,$ the following
conditions hold:

\begin{itemize}
\item[(b$_{1}$)] $d(x,y)=0$ if and only if $x=y,$

\item[(b$_{2}$)] $d(x,y)=d(y,x),$

\item[(b$_{3}$)] $d(x,y)\leq b\left( d(x,z)+d(z,y)\right) ,$
\end{itemize}

\noindent The pair $(X,d)$ is called a $b$- metric space with parameter $%
b\geq 1$.

\noindent If $b=1,$ then $b$-metric space is a metric spaces. But the
converse does not hold in general \cite{AN, C4, CZ}.\medskip

\noindent \textbf{Example 1.2. }\cite{ROS} Let $(X,d)$ be a metric space,
and $\rho (x,y)=(d(x,y))^{p},$ where $p>1$ is a real number. Then $\rho $ is
a $b$-metric with $b=2^{p-1}.$

\noindent Obviously conditions (b$_{1}$) and (b$_{2}$)\ of above definition
are satisfied. If $1<p<\infty ,$ then the convexity of the function $%
f(x)=x^{p}$ $(x>0)$ implies%
\begin{equation*}
\left( \frac{a+b}{2}\right) ^{p}\leq \frac{1}{2}\left( a^{p}+b^{p}\right) ,
\end{equation*}%
and hence, $\left( a+b\right) ^{p}\leq 2^{p-1}(a^{p}+b^{p})$ holds. Thus,
for each $x,y,z\in X$ we obtain%
\begin{eqnarray*}
\rho (x,y) &=&(d(x,y))^{p}\leq (d(x,z)+d(z,y))^{p} \\
&\leq &2^{p-1}\left( (d(x,z))^{p}+(d(z,y))^{p}\right) \\
&=&2^{p-1}(\rho (x,z)+\rho (z,y)).
\end{eqnarray*}%
So condition (b$_{3}$) of the above definition is satisfied and $\rho $ is a 
$b$-metric.

\noindent If $X=%
\mathbb{R}
$ ( set of real numbers ) and $d(x,y)=\left\vert x-y\right\vert $ is the
usual metric, then $\rho (x,y)=(x-y)^{2}$ is a $b$-metric on $%
\mathbb{R}
$ with $b=2,$ but is not a metric on $%
\mathbb{R}
$.\medskip

\noindent \textbf{Definition 1.3.} \ \ A set $C$ is compact if every
sequence $\{x_{n}\}$ in $C$ contains a subsequence having a limit in $C$.
Note that closed and bounded subsets of $%
\mathbb{R}
^{n}$ are compact. Also, every finite set in $%
\mathbb{R}
^{n}$ are compact. On the other hands, $(0,1]\subset 
\mathbb{R}
$ is not compact as $\{1,\dfrac{1}{2},\dfrac{1}{2^{2}},...\}\subset (0,1]$
does not have any convergent subsequence. Also, $Z\subset 
\mathbb{R}
$ is not compact.

Let $\mathcal{H}(X)$ denotes the set of all non-empty compact subsets of $X$%
. For $A,B\in \mathcal{H}(X)$, let%
\begin{equation*}
H(A,B)=\max \{\sup\limits_{b\in B}d(b,A),\sup\limits_{a\in A}d(a,B)\},
\end{equation*}%
where $d(x,B)=\inf \{d(x,y):y\in B\}$ is the distance of a point $x$ from
the set $B$. The mapping $H$ is said to be the Pompeiu-Hausdorff metric
induced by $d$. If $(X,d)$ is a compete $b$-metric space, then $\left( 
\mathcal{H}(X),H\right) $ is also a complete $b$-metric space.\medskip

For the sake of completeness, we state the following Lemma hold in $b$%
-metric space \cite{NSA16}.

\noindent \textbf{Lemma 1.4.} \ \ Let $(X,d)$ be a $b$-metric space. For all 
$A,B,C,D\in \mathcal{H}(X)$, the following hold:

\begin{description}
\item[(i)] If $B\subseteq C,$ then $\sup\limits_{a\in A}d(a,C)\leq
\sup\limits_{a\in A}d(a,B).$

\item[(ii)] $\sup\limits_{x\in A\cup B}d(x,C)=\max \{\sup\limits_{a\in
A}d(a,C),\sup\limits_{b\in B}d(b,C)\}.$

\item[(iii)] $H(A\cup B,C\cup D)\leq \max \{H(A,C),H(B,D)\}.$
\end{description}

\noindent The following lemmas from \cite{CZ, CZ1, CZ2} will be needed in
the sequel to prove the main result of the paper.

\noindent \textbf{Lemma 1.5.} \ Let $(X,d)$ be a $b$-metric space$,$ $x,y\in
X$ and $A,B\in CB(X).$ The following statements hold:

\begin{enumerate}
\item $(CB(X),H)$ is a $b$-metric space$.$

\item $d(x,B)\leq H(A,B)$ for all $x\in A.$

\item $d(x,A)\leq b\left( d(x,y)+d(y,A)\right) .$

\item For $h>1$ and $\acute{a}\in A,$ there is a $\acute{b}\in B$ such that $%
d(\acute{a},\acute{b})\leq hH(A,B).$

\item For every $h>0$ and $\acute{a}\in A,$ there is a $\acute{b}\in B$ such
that $d(\acute{a},\acute{b})\leq H(A,B)+h.$

\item For every $\lambda >0$ and $\tilde{a}\in A,$ there is a $\tilde{b}\in
B $ such that $d(\tilde{a},\tilde{b})\leq \lambda .$

\item For every $\lambda >0$ and $\tilde{a}\in A,$ there is a $\tilde{b}\in
B $ such that $d(\tilde{a},\tilde{b})\leq \lambda \ $implies$\ H(A,B)\leq
\lambda .$

\item $d(x,A)=0$ if and only if $x\in \bar{A}=A.$

\item For $\{x_{n}\}\subseteq X$,%
\begin{equation*}
d(x_{0},x_{n})\leq
bd(x_{0},x_{1})+...+b^{n-1}d(x_{n-2},x_{n-1})+b^{n-1}d(x_{n-1},x_{n}).
\end{equation*}
\end{enumerate}

\noindent \textbf{Definition 1.6.} \ \ Let $(X,d)$ be a $b$-metric space. A
sequence $\{x_{n}\}$ in $X$ is called:

\begin{description}
\item[(i)] Cauchy if and only if for $\varepsilon >0,$ there exists $%
n(\varepsilon )\in 
\mathbb{N}
$ such that for each $n,m\geq n(\varepsilon )$ we have $d(x_{n},x_{m})<%
\varepsilon .$

\item[(ii)] Convergent if and only if there exists $x\in X$ such that for
all $\varepsilon >0$ there exists $n(\varepsilon )\in 
\mathbb{N}
$ such that for all $n\geq n(\varepsilon )$ we have $d(x_{n},x)<\varepsilon $%
. In this case we write $\lim\limits_{n\rightarrow \infty }x_{n}=x.$
\end{description}

It is known that a sequence $\{x_{n}\}$ in $b$-metric space $X$ is Cauchy if
and only if $\lim\limits_{n\rightarrow \infty }d(x_{n},x_{n+p})=0$ for all $%
p\in 
\mathbb{N}
$. A sequence $\{x_{n}\}$ is convergent to $x\in X$ if and only if $%
\lim\limits_{n\rightarrow \infty }d(x_{n},x)=0.$ A subset $Y\subset X$ is
closed if and only if for each sequence $\{x_{n}\}$ in $Y$ that converges to
an element $x$, we have $x\in Y.$ A subset $Y\subset X$ is bounded if $%
diam(Y)$ is finite, where $diam(Y)=\sup \left\{ d(a,b),a,b\in Y\right\} $. A 
$b$-metric space $(X,d)$ is said to be complete if every Cauchy sequence in $%
X$ is convergent in $X$.\medskip

An et al. \cite{AN} studied the topological properties of $b$-metric spaces
and stated the following assertions:

\begin{description}
\item[(c$_{1}$)] In a $b$-metric space $(X,d),$ $d$ is not necessarily
continuous in each variable.

\item[(c$_{2}$)] In a $b$-metric space $(X,d),$ If $d$ is continuous in one
variable then $d$ is continuous in other variable.

\item[(c$_{3}$)] An open ball in $b$-metric space $(X,d)$ is not necessarily
an open set. An open ball is open if $d$ is continuous in one variable.
\end{description}

Wardowski \cite{Wardowski12} introduced another generalized contraction
called $F$-contraction and proved a fixed point result as an interesting
generalization of the Banach contraction principle in complete metric space
(See also \cite{KW15}).

Let $\digamma $ be the collection of all continuous mappings $F:%
\mathbb{R}
^{+}\rightarrow 
\mathbb{R}
$ that satisfy the following conditions:

\begin{itemize}
\item[($F_{1}$)] $F$ is strictly increasing, that is, for all $\alpha ,\beta
\in 
\mathbb{R}
^{+}$\ such that $\alpha <\beta $ implies that $F(\alpha )<F(\beta )$.

\item[($F_{2}$)] For every sequence $\{ \alpha _{n}\}$ of positive real
numbers, $\lim \limits_{n\rightarrow \infty }\alpha _{n}=0$ and $\lim
\limits_{n\rightarrow \infty }F\left( \alpha _{n}\right) =-\infty $ are
equivalent.

\item[($F_{3}$)] There exists $h\in \left( 0,1\right) $ such that $%
\lim\limits_{\alpha \rightarrow 0^{+}}\alpha ^{h}F(\alpha )=0$.
\end{itemize}

\noindent \textbf{Definition 1.7.} \cite{Wardowski12} Let $(X,d)$ be a
metric space. A self-mapping $f$ on $X$ is called an $F$-contraction if for
any $x,y\in X$, there exists $F\in \digamma $ and $\tau >0$ such that%
\begin{equation}
\tau +F(d(fx,fy))\leq F(d(x,y)),  \tag{1.1}
\end{equation}%
whenever $d(fx,fy)>0$.

\noindent From ($F_{1}$) and (1.1), we conclude that%
\begin{equation*}
d(fx,fy)<d(x,y),\text{ for all }x,y\in X,\text{ }fx\neq fy,
\end{equation*}%
that is, every $F$-contraction mapping is contractive, and in particular,
every $F$-contraction mapping is continuous.\medskip

Wardowski \cite{Wardowski12} proved that in complete metric space $(X,d)$,
every $F$-contractive self-map has a unique fixed point in $X$ and for every 
$x_{0}$ in $X$ a sequence of iterates $\{x_{0},fx_{0},f^{2}x_{0},...\}$
converges to the fixed point of $f$.\medskip

Let $\Upsilon $ be the set of all mapping $\tau :%
\mathbb{R}
_{+}\rightarrow 
\mathbb{R}
_{+}$ that satisfying $\lim \inf_{t\rightarrow 0}\tau (t)>0$ for all $t\geq
0.$\medskip

\noindent \textbf{Definition 1.8. \ }\ Let $(X,d)$ be a $b$-metric space. A
self-mapping $f$ on $X$ is called a generalized $F$-contraction if for any $%
x,y\in X$, there exists $F\in \digamma $ and $\tau \in \Upsilon $ such that%
\begin{equation}
\tau (d\left( x,y\right) )+F(d(fx,fy))\leq F(d(x,y)),  \tag{1.2}
\end{equation}%
whenever $d(fx,fy)>0$.\medskip

\noindent \textbf{Theorem 1.9.} \ \ Let $(X,d)$ be a $b$-metric space and $%
f:X\rightarrow X$ an generalized $F$-contraction. Then

\begin{description}
\item[(1)] $f$ maps elements in $\mathcal{H}(X)$ to elements in $\mathcal{H}%
(X)$.

\item[(2)] If for any $A\in \mathcal{H}(X),$%
\begin{equation*}
f(A)=\{f(x):x\in A\}.\text{ }
\end{equation*}%
Then $f:\mathcal{H}(X)\rightarrow \mathcal{H}(X)$ is a generalized $F$%
-contraction mapping on $(\mathcal{H}(X),H)$.
\end{description}

\noindent \textit{Proof}. \ \ As generalized $F$-contractive mapping is
continuous and the image of a compact subset under $f:X\rightarrow X$ \ is
compact, so we obtain%
\begin{equation*}
A\in \mathcal{H}(X)\text{ implies }f(A)\in \mathcal{H}(X).
\end{equation*}%
To prove (2): Let $A,B\in \mathcal{H}(X)$ with $H\left( f\left( A\right)
,f\left( B\right) \right) \neq \emptyset $. Since $f:X\rightarrow X$ is a
generalized $F$-contraction, we obtain%
\begin{equation*}
0<d\left( fx,fy\right) <d\left( x,y\right) \text{ for all }x,y\in X,\text{ }%
x\neq y.
\end{equation*}%
Thus we have%
\begin{equation*}
d\left( fx,f\left( B\right) \right) =\inf_{y\in B}d\left( fx,fy\right)
<\inf_{y\in B}d\left( x,y\right) =d\left( x,B\right) .
\end{equation*}%
Also%
\begin{equation*}
d\left( fy,f\left( A\right) \right) =\inf_{x\in A}d\left( fy,fx\right)
<\inf_{x\in A}d\left( y,x\right) =d\left( y,A\right) .
\end{equation*}%
Now%
\begin{eqnarray*}
H\left( f\left( A\right) ,f\left( B\right) \right) &=&\max
\{\sup\limits_{x\in A}d(fx,f\left( B\right) ),\sup\limits_{y\in
B}d(fy,f\left( A\right) )\} \\
&<&\max \{\sup\limits_{x\in A}d(x,B),\sup\limits_{y\in B}d(y,A)\}=H\left(
A,B\right) .
\end{eqnarray*}%
By strictly increasing of $F$ implies%
\begin{equation*}
F\left( H\left( f\left( A\right) ,f\left( B\right) \right) \right) <F\left(
H\left( A,B\right) \right) .
\end{equation*}%
Consequently, there exists a function $\tau :%
\mathbb{R}
_{+}\rightarrow 
\mathbb{R}
_{+}$ with $\lim \inf_{t\rightarrow 0}\tau (t)>0$ for all $t\geq 0$ such that%
\begin{equation*}
\tau \left( H\left( A,B\right) \right) +F\left( H\left( f\left( A\right)
,f\left( B\right) \right) \right) \leq F\left( H\left( A,B\right) \right) .
\end{equation*}%
Hence $f:\mathcal{H}(X)\rightarrow \mathcal{H}(X)$ is a generalized $F$%
-contraction. $\square $\medskip

\noindent \textbf{Theorem 1.10. \ } Let $(X,d)$ be a $b$-metric space and $%
\{f_{n}:n=1,2,...,N\}$ a finite family of generalized $F$-contraction
self-mappings on $X.$ Define $T:\mathcal{H}(X)\rightarrow \mathcal{H}(X)$ by%
\begin{eqnarray*}
T(A) &=&f_{1}(A)\cup f_{2}(A)\cup \cdot \cdot \cdot \cup f_{N}(A) \\
&=&\cup _{n=1}^{N}f_{n}(A),\text{ for each }A\in \mathcal{H}(X).
\end{eqnarray*}%
Then $T$ is a generalized $F$-contraction on $\mathcal{H}\left( X\right) $.

\noindent \textit{Proof.} \ \ \ We demonstrate the claim for $N=2$. Let $%
f_{1},f_{2}:X\rightarrow X$ be two $F$-contractions. Take $A,B\in \mathcal{H}%
\left( X\right) $ with $H(T(A),T(B))\neq 0.$ From Lemma 1.4 (iii), it
follows that%
\begin{eqnarray*}
&&\tau \left( H\left( A,B\right) \right) +F\left( H(T(A),T(B))\right) \\
&=&\tau \left( H\left( A,B\right) \right) +F\left( H(f_{1}(A)\cup
f_{2}(A),f_{1}(B)\cup f_{2}(B))\right) \\
&\leq &\tau \left( H\left( A,B\right) \right) +F\left( \max
\{H(f_{1}(A),f_{1}(B)),H(f_{2}(A),f_{2}(B))\}\right) \\
&\leq &F\left( H(A,B)\right) .\text{ \ }\square
\end{eqnarray*}

\noindent \textbf{Definition 1.11.} \ \ Let $(X,d)$ be a metric space. A
mapping $T:\mathcal{H}\left( X\right) \rightarrow \mathcal{H}\left( X\right) 
$ is said to be a Ciric type generalized $F$-contraction if for $F\in
\digamma $ and $\tau \in \Upsilon $ such that for any $A$, $B\in \mathcal{H}%
\left( X\right) $ with $H(T(A),T(B))\neq 0$, the following holds:%
\begin{equation}
\tau \left( M_{T}\left( A,B\right) \right) +F\left( H\left( T\left( A\right)
,T\left( B\right) \right) \right) \leq F(M_{T}(A,B)),  \tag{1.3}
\end{equation}%
where%
\begin{eqnarray*}
M_{T}(A,B) &=&\max \{H(A,B),H(A,T\left( A\right) ),H(B,T\left( B\right) ),%
\dfrac{H(A,T\left( B\right) )+H(B,T\left( A\right) )}{2b}, \\
&&H(T^{2}\left( A\right) ,T\left( A\right) ),H(T^{2}\left( A\right)
,B),H(T^{2}\left( A\right) ,T\left( B\right) )\}.
\end{eqnarray*}

\noindent \textbf{Theorem 1.12. \ } Let $(X,d)$ be a $b$-metric space and $%
\{f_{n}:n=1,2,...,N\}$ a finite sequence of generalized $F$-contraction
mappings on $X.$ If $T:\mathcal{H}(X)\rightarrow \mathcal{H}(X)$ is defined
by%
\begin{eqnarray*}
T(A) &=&f_{1}(A)\cup f_{2}(A)\cup \cdot \cdot \cdot \cup f_{N}(A) \\
&=&\cup _{n=1}^{N}f_{n}(A),\text{ for each }A\in \mathcal{H}(X).
\end{eqnarray*}%
Then $T$ is a Ciric type generalized $F$-contraction mapping on $\mathcal{H}%
\left( X\right) $.

\noindent \textit{Proof.} \ \ \ Using Theorem 1.11 with property $(F_{1})$,
the result follows. $\square $ \medskip

An operator $T$ in above Theorem is called Ciric type generalized $F$%
-Hutchinson operator.

\noindent \textbf{Definition 1.13.} \ \ Let $X$ be a complete $b$-metric
space. If $f_{n}:X\rightarrow X$, $n=1,2,...,N$ are generalized $F$%
-contraction mappings, then $(X;f_{1},f_{2},...,f_{N})$ is called
generalized $F$-contractive iterated function system (IFS).\medskip

Thus generalized $F$-contractive iterated function system consists of a
complete $b$-metric space and finite family of generalized $F$-contraction
mappings on $X.$

\noindent \textbf{Definition 1.14.} \ \ A nonempty compact set $A\subset X$
is said to be an attractor of the generalized $F$-contractive IFS $T$ if

\begin{itemize}
\item[(a)] $T(A)=A$ and

\item[(b)] there is an open set $U\subset X$ such that $A\subset U$ and $%
\lim\limits_{n\rightarrow \infty }T^{n}(B)=A$ for any compact set $B\subset
U $, where the limit is taken with respect to the Hausdorff metric.
\end{itemize}

\section{Main Results}

\noindent We start with the following result.

\noindent \textbf{Theorem 2.1. \ } Let $(X,d)$ be a complete $b$-metric
space and $\{X;f_{n},n=1,2,$\textperiodcentered \textperiodcentered
\textperiodcentered $,k\}$ a generalized $F$-contractive iterated function
system.$\ $Then following hold:

\noindent (a) $\ $A mapping $T:\mathcal{H}(X)\rightarrow \mathcal{H}(X)$
defined by%
\begin{equation*}
T(A)=\cup _{n=1}^{k}f_{n}(A),\text{ for all }A\in \mathcal{H}(X)
\end{equation*}%
is Ciric type generalized $F$-Hutchinson operator.

\noindent (b) \ Operator $T$ has a unique fixed point $U\in \mathcal{H}%
\left( X\right) ,$ that is%
\begin{equation*}
U=T\left( U\right) =\cup _{n=1}^{k}f_{n}(U).
\end{equation*}%
(c) \ For any initial set $A_{0}\in \mathcal{H}\left( X\right) $, the
sequence of compact sets $\{A_{0},T\left( A_{0}\right) ,T^{2}\left(
A_{0}\right) ,...\}$ converges to a fixed point of $T$.

\noindent \textit{Proof. \ \ }\ Part (a) follows from Theorem 1.13. For
parts (b) and (c), we proceed as follows: \ \textit{\ }Let $A_{0}\ $be an
arbitrary element in $\mathcal{H}\left( X\right) .$ If $A_{0}=T\left(
A_{0}\right) ,$ then the proof is finished. So we assume that $A_{0}\neq
T\left( A_{0}\right) .$ Define%
\begin{equation*}
A_{1}=T(A_{0}),\text{ }A_{2}=T\left( A_{1}\right) ,...,A_{m+1}=T\left(
A_{m}\right)
\end{equation*}%
for $m\in 
\mathbb{N}
.$

We may assume that $A_{m}\neq A_{m+1}$ for all $m\in 
\mathbb{N}
.$ If not, then $A_{k}=A_{k+1}$ for some $k$ implies $A_{k}=T(A_{k})\ $and
this completes the proof. Take $A_{m}\neq A_{m+1}$ for all $m\in 
\mathbb{N}
$. From (1.7), we have%
\begin{eqnarray*}
&&\tau \left( M_{T}\left( A_{m},A_{m+1}\right) \right) +F\left(
H(A_{m+1},A_{m+2})\right) \\
&=&\tau (M_{T}\left( A_{m},A_{m+1}\right) )+F\left( H(T\left( A_{m}\right)
,T\left( A_{m+1}\right) )\right) \\
&\leq &F\left( M_{T}\left( A_{m},A_{m+1}\right) \right) ,
\end{eqnarray*}%
where%
\begin{eqnarray*}
M_{T}\left( A_{m},A_{m+1}\right) &=&\max \{H(A_{m},A_{m+1}),H\left(
A_{m},T\left( A_{m}\right) \right) ,H\left( A_{m+1},T\left( A_{m+1}\right)
\right) , \\
&&\frac{H\left( A_{m},T\left( A_{m+1}\right) \right) +H\left(
A_{m+1},T\left( A_{m}\right) \right) }{2b}, \\
&&H(T^{2}\left( A_{m}\right) ,T\left( A_{m}\right) ),H\left( T^{2}\left(
A_{m}\right) ,A_{m+1}\right) ,H\left( T^{2}\left( A_{m}\right) ,T\left(
A_{m+1}\right) \right) \} \\
&=&\max \{H(A_{m},A_{m+1}),H\left( A_{m},A_{m+1}\right) ,H\left(
A_{m+1},A_{m+2}\right) , \\
&&\frac{H\left( A_{m},A_{m+2}\right) +H\left( A_{m+1},A_{m+1}\right) }{2b},
\\
&&H(A_{m+2},A_{m+1}),H\left( A_{m+2},A_{m+1}\right) ,H\left(
A_{m+2},A_{m+2}\right) \} \\
&=&\max \{H(A_{m},A_{m+1}),H\left( A_{m+1},A_{m+2}\right) \}.
\end{eqnarray*}%
In case $M_{T}\left( A_{m},A_{m+1}\right) =H(A_{m+1},A_{m+2}),$ we have%
\begin{equation*}
F\left( H(A_{m+1},A_{m+2})\right) \leq F(H(A_{m+1},A_{m+2}))-\tau
(H(A_{m+1},A_{m+2})),
\end{equation*}%
a contradiction as $\tau (H(A_{m+1},A_{m+2}))>0$. Therefore $M_{T}\left(
A_{m},A_{m+1}\right) =H(A_{m},A_{m+1})$ and we have%
\begin{eqnarray*}
F\left( H(A_{m+1},A_{m+2})\right) &\leq &F(H(A_{m},A_{m+1}))-\tau
(H(A_{m},A_{m+1})) \\
&<&F(H(A_{m},A_{m+1})).
\end{eqnarray*}%
Thus $\{H(A_{m+1},A_{m+2})\}$ is decreasing and hence convergent. We now
show that $\lim\limits_{m\rightarrow \infty }H(A_{m+1},A_{m+2})=0.$ By
property of $\tau ,$ there exists $c>0$ with $n_{0}\in 
\mathbb{N}
$ such that $\tau (H(A_{m},A_{m+1}))$ $>c$ for all $m\geq n_{0}.$ Note that%
\begin{eqnarray*}
F(H(A_{m+1},A_{m+2})) &\leq &F\left( H(A_{m},A_{m+1})\right) -\tau \left(
H(A_{m},A_{m+1})\right) \\
&\leq &F(H(A_{m-1},A_{m}))-\tau (H(A_{m-1},A_{m}))-\tau \left(
H(A_{m},A_{m+1})\right) \\
&\leq &...\leq H(A_{0},A_{1})-[\tau \left( H(A_{0},A_{1})\right) +\tau
(H(A_{1},A_{2})) \\
&&+...+\tau (H(A_{m},A_{m+1})] \\
&\leq &F(H(A_{0},A_{1}))-n_{0},
\end{eqnarray*}%
gives $\lim\limits_{m\rightarrow \infty }F\left( H(A_{m+1},A_{m+2})\right)
=-\infty $ which together with ($F_{2}$) implies that $\lim\limits_{m%
\rightarrow \infty }H(A_{m+1},A_{m+2})=0.$ By ($F_{3}$), there exists $h\in
\left( 0,1\right) $ such that%
\begin{equation*}
\lim\limits_{n\rightarrow \infty }[H(A_{m+1},A_{m+2})]^{h}F\left(
H(A_{m+1},A_{m+2})\right) =0.
\end{equation*}%
Thus we have%
\begin{eqnarray*}
&&[H(A_{m},A_{m+1})]^{h}F\left( H(A_{m},A_{m+1})\right)
-[H(A_{m},A_{m+1})]^{h}F\left( H(A_{0},A_{1})\right) \\
&\leq
&[H(A_{m},A_{m+1})]^{h}(F(H(A_{0},A_{1})-n_{0}))-[H(A_{m},A_{m+1})]^{h}F%
\left( H(A_{0},A_{1})\right) \\
&\leq &-n_{0}[H(A_{m},A_{m+1})]^{h}\leq 0.
\end{eqnarray*}%
On taking limit as $n\rightarrow \infty $ we obtain that $%
\lim\limits_{m\rightarrow \infty }m[H(A_{m+1},A_{m+2})]^{h}=0.$ Hence $%
\lim\limits_{m\rightarrow \infty }m^{\frac{1}{h}}H(A_{m+1},A_{m+2})=0.$
There exists $n_{1}\in 
\mathbb{N}
$ such that $m^{\frac{1}{h}}H(A_{m+1},A_{m+2})\leq 1$ for all $m\geq n_{1}$
and hence $H(A_{m+1},A_{m+2})\leq \frac{1}{m^{1/h}}$ for all $m\geq n_{1}.$
For $m,n\in 
\mathbb{N}
$ with $m>n\geq n_{1}$, we have%
\begin{eqnarray*}
H\left( A_{n},A_{m}\right) &\leq &H\left( A_{n},A_{n+1}\right) +H\left(
A_{n+1},A_{n+2}\right) +...+H\left( A_{m-1},A_{m}\right) \\
&\leq &\sum_{i=n}^{\infty }\frac{1}{i^{1/h}}.
\end{eqnarray*}%
By the convergence of the series $\sum_{i=1}^{\infty }\dfrac{1}{i^{1/h}},$
we get $H\left( A_{n},A_{m}\right) \rightarrow 0$ as $n,m\rightarrow \infty $%
. Therefore $\{A_{n}\}$ is a Cauchy sequence in $X.$ Since $(\mathcal{H}%
(X),d)$ is complete, we have $A_{n}\rightarrow U$ as $n\rightarrow \infty $\
for some $U\in \mathcal{H}(X).$

In order to show that $U$ is the fixed point of $T,$ we contrary assume that
Pompeiu-Hausdorff weight assign to the $U$ and $T\left( U\right) $ is not
zero. Now%
\begin{eqnarray}
&&\tau \left( M_{T}\left( A_{n},U\right) \right) +F\left( H(A_{n+1},T\left(
U\right) )\right)  \notag \\
&=&\tau +F(H(T\left( A_{n}\right) ,T\left( U\right) ))\leq F\left(
M_{T}\left( A_{n},U\right) \right) ,  \TCItag{2.1}
\end{eqnarray}%
where%
\begin{eqnarray*}
M_{T}\left( A_{n},U\right) &=&\max \{H(A_{n},U),H(A_{n},T\left( A_{n}\right)
),H(U,T\left( U\right) ),\frac{H(A_{n},T\left( U\right) ){\small +}%
H(U,T\left( A_{n}\right) )}{2b}, \\
&&H(T^{2}\left( A_{n}\right) ,T\left( A_{n}\right) ),H(T^{2}\left(
A_{n}\right) ,U),H(T^{2}\left( A_{n}\right) ,T\left( U\right) )\} \\
&=&\max \{H(A_{n},U),H(A_{n},A_{n+1}),H(U,T\left( U\right) ),\frac{%
H(A_{n},T\left( U\right) ){\small +}H(U,A_{n+1})}{2b}, \\
&&H(A_{n+2},A_{n+1}),H(A_{n+2},U),H(A_{n+2},T\left( U\right) )\}\}.
\end{eqnarray*}%
Now we consider the following cases:

\noindent (i). If $M_{T}\left( A_{n},U\right) =H(A_{n},U),$ then on taking
lower limit as $n\rightarrow \infty $ in (2.1), we have%
\begin{equation*}
\underset{n\rightarrow \infty }{\lim \inf }\tau (H(A_{n},U))+F\left(
H(T\left( U\right) ,U)\right) \leq F\left( H\left( U,U\right) \right) ,
\end{equation*}%
a contradiction as $\lim \inf_{t\rightarrow 0}\tau (t)>0$ for all $t\geq 0$.

\noindent (2). When $M_{T}\left( A_{n},U\right) =H(A_{n},A_{n+1}),$ then by
taking lower limit as $n\rightarrow \infty ,$ we obtain%
\begin{equation*}
\underset{n\rightarrow \infty }{\lim \inf }\tau (H(A_{n},A_{n+1}))+F\left(
H(T\left( U\right) ,U)\right) \leq F\left( H\left( U,U\right) \right) ,
\end{equation*}%
gives a contradiction.

\noindent (3). In case $M_{T}\left( A_{n},U\right) =H(U,T\left( U\right) ),$
then we get%
\begin{equation*}
\tau (H(U,T\left( U\right) ))+F\left( H(T\left( U\right) ,U)\right) \leq
F\left( H\left( U,T\left( U\right) \right) \right) ,
\end{equation*}%
a contradiction as $\tau (H(U,T\left( U\right) ))>0$.

\noindent (4). If $M_{T}\left( A_{n},U\right) =\dfrac{H(A_{n},T\left(
U\right) )+H(U,A_{n+1})}{2b},$ then on taking lower limit as $n\rightarrow
\infty ,$ we have%
\begin{eqnarray*}
&&\underset{n\rightarrow \infty }{\lim \inf }\tau (\dfrac{H(A_{n},T\left(
U\right) )+H(U,A_{n+1})}{2b})+F\left( H(T\left( U\right) ,U)\right) \\
&\leq &F(\frac{H\left( U,T\left( U\right) \right) +H\left( U,U\right) }{2b})
\\
&=&F(\frac{H\left( U,T\left( U\right) \right) }{2b}),
\end{eqnarray*}%
a contradiction as $F$ is strictly increasing map.

\noindent (5). When $M_{T}\left( A_{n},U\right) =H(A_{n+2},A_{n+1}),$ then%
\begin{equation*}
\underset{n\rightarrow \infty }{\lim \inf }\tau (H(A_{n+2},A_{n+1}))+F\left(
H(T\left( U\right) ,U)\right) \leq F\left( H\left( U,U\right) \right) ,
\end{equation*}%
gives a contradiction.

\noindent (6). In case $M_{T}\left( A_{n},U\right) =H(A_{n+2},U),$ then on
taking lower limit as $n\rightarrow \infty $ in (2.1), we get%
\begin{equation*}
\underset{n\rightarrow \infty }{\lim \inf }\tau (H(A_{n+2},U))+F\left(
H(T\left( U\right) ,U)\right) \leq F\left( H\left( U,U\right) \right) ,
\end{equation*}%
a contradiction.

\noindent (7). Finally if $M_{T}\left( A_{n},U\right) =H(A_{n+2},T\left(
U\right) ),$ then on taking lower limit as $n\rightarrow \infty ,$ we have%
\begin{equation*}
\underset{n\rightarrow \infty }{\lim \inf }\tau (H(A_{n+2},T\left( U\right)
))+F\left( H(T\left( U\right) ,U)\right) \leq F(H(U,T\left( U\right) )),
\end{equation*}%
a contradiction.

\noindent Thus, $U$ is the fixed point of $T$.

\noindent To show the uniqueness of fixed point of $T$, assume that $U$ and $%
V$ are two fixed points of $T$ with $H\left( U,V\right) $ is not zero. Since 
$T$ is a $F$-contraction map, we obtain that%
\begin{eqnarray*}
\tau \left( M_{T}\left( U,V\right) \right) +F(H(U,V)) &=&\tau \left(
M_{T}\left( U,V\right) \right) +F(H(T\left( U\right) ,T\left( V\right) )) \\
&\leq &F\left( M_{T}\left( U,V\right) \right) ,
\end{eqnarray*}%
where%
\begin{eqnarray*}
M_{T}\left( U,V\right) &=&\max \{H(U,V),H(U,T\left( U\right) ),H(V,T\left(
V\right) ),\frac{H(U,T\left( V\right) )+H(V,T\left( U\right) )}{2b}, \\
&&H(T^{2}\left( U\right) ,U),H(T^{2}\left( U\right) ,V),H(T^{2}\left(
U\right) ,T\left( V\right) )\} \\
&=&\max \{H\left( U,V\right) ,H(U,U),H(V,V),\frac{H(U,V)+H(V,U)}{2b}, \\
&&H\left( U,U\right) ,H(U,V),H(U,V)\} \\
&=&H(U,V),
\end{eqnarray*}%
that is,%
\begin{equation*}
\tau \left( H\left( U,V\right) \right) +F(H(U,V))\leq F\left( H\left(
U,V\right) \right) ,
\end{equation*}%
a contradiction as $\tau \left( H\left( U,V\right) \right) >0$. Thus $T$ has
a unique fixed point $U\in \mathcal{H}(X)$. $\square $\medskip

\noindent \textbf{Remark 2.2.} \ \ In Theorem 2.1, if we take $\mathcal{S}%
(X) $ the collection of all singleton subsets of $X,$ then clearly $\mathcal{%
S}(X)\subseteq \mathcal{H}(X).$ Moreover, consider $f_{n}=f$ for each $n,$
where $f=f_{i}$ for any $i\in \{1,2,3,...,k\},$ then the mapping $T$ becomes%
\begin{equation*}
T(x)=f(x).
\end{equation*}%
With this setting we obtain the following fixed point result. \medskip

\noindent \textbf{Corollary 2.3.} \ \ Let $(X,d)$ be a complete $b$-metric
space and $\{X:f_{n},n=1,2,$\textperiodcentered \textperiodcentered
\textperiodcentered $,k\}$ a generalized iterated function system.\ Let $%
f:X\rightarrow X$ be a mapping defined as in Remark 2.2. If there exist some 
$F\in \digamma $ and $\tau \in \Upsilon $ such that for any $x,y\in \mathcal{%
H}\left( X\right) $ with $d(f(x),f(y))\neq 0$, the following holds:%
\begin{equation*}
\tau \left( M_{f}(x,y)\right) +F\left( d\left( fx,fy\right) \right) \leq
F(M_{f}(x,y)),
\end{equation*}%
where%
\begin{eqnarray*}
M_{T}(x,y) &=&\max \{d(x,y),d(x,fx),d(y,fy),\dfrac{d(x,fy)+d(y,fx)}{2b}, \\
&&d(f^{2}x,y),d(f^{2}x,fx),d(f^{2}x,fy)\}.
\end{eqnarray*}%
Then $f$ has a unique fixed point in $X.$ Moreover, for any initial set $%
x_{0}\in X$, the sequence of compact sets $\{x_{0},fx_{0},f^{2}x_{0},...\}$
converges to a fixed point of $f$. \medskip

\noindent \textbf{Corollary 2.4. \ } Let $(X,d)$ be a complete $b$-metric
space and $(X;f_{n},n=1,2,$\textperiodcentered \textperiodcentered
\textperiodcentered $,k)$ be iterated function system where each $f_{i}$ for 
$i=1,2,...,k$ is a contraction self-mapping on $X.\ $Then $T:\mathcal{H}%
(X)\rightarrow \mathcal{H}(X)$ defined in Theorem 2.1 has a unique fixed
point in $\mathcal{H}\left( X\right) .$ Furthermore, for any set $A_{0}\in 
\mathcal{H}\left( X\right) $, the sequence of compact sets $\{A_{0},T\left(
A_{0}\right) ,T^{2}\left( A_{0}\right) ,...\}$ converges to a fixed point of 
$T$.

\noindent \textit{Proof.} \ \ \ It follows from Theorem 1.10 that if each $%
f_{i}$ for $i=1,2,...,k$ is a contraction mapping on $X,$ then the mapping $%
T:\mathcal{H}(X)\rightarrow \mathcal{H}(X)$ defined by%
\begin{equation*}
T(A)=\cup _{n=1}^{k}f_{n}(A),\text{ for all }A\in \mathcal{H}(X)
\end{equation*}%
is contraction on $\mathcal{H}\left( X\right) $. Using Theorem 2.1, the
result follows. $\square $\medskip

\noindent \textbf{Corollary 2.5. \ } Let $(X,d)$ be a complete $b$-metric
space and $(X;f_{n},n=1,2,$\textperiodcentered \textperiodcentered
\textperiodcentered $,k)$ an iterated function system. Suppose that each $%
f_{i}$ for $i=1,2,...,k$ is a mapping on $X$ satisfying%
\begin{equation*}
d\left( f_{i}x,f_{i}y\right) e^{d\left( f_{i}x,f_{i}y\right) -d\left(
x,y\right) }\leq e^{-\tau \left( d(x,y)\right) }d\left( x,y\right) ,
\end{equation*}%
for all $x,y\in X,$ $f_{i}x\neq f_{i}y,$ where $\tau \in \Upsilon .$ Then
the mapping $T:\mathcal{H}(X)\rightarrow \mathcal{H}(X)$ defined in Theorem
2.1 has a unique fixed point in $\mathcal{H}\left( X\right) .\ $Furthermore,
for any set $A_{0}\in \mathcal{H}\left( X\right) $, the sequence of compact
sets $\{A_{0},T\left( A_{0}\right) ,T^{2}\left( A_{0}\right) ,...\}$
converges to a fixed point of $T$.

\noindent \textit{Proof.} \ \ \ Take $F\left( \lambda \right) =\ln \left(
\lambda \right) +\lambda ,$ $\lambda >0$ in Theorem 1.10, then each mapping $%
f_{i}$ for $i=1,2,...,k$ on $X$ satisfies%
\begin{equation*}
d\left( f_{i}x,f_{i}y\right) e^{d\left( f_{i}x,f_{i}y\right) -d\left(
x,y\right) }\leq e^{-\tau \left( d(x,y)\right) }d\left( x,y\right) ,
\end{equation*}%
for all $x,y\in X,$ $f_{i}x\neq f_{i}y,$ where $\tau \in \Upsilon .$ Again
form Theorem 1.10, the mapping $T:\mathcal{H}(X)\rightarrow \mathcal{H}(X)$
defined by%
\begin{equation*}
T(A)=\cup _{n=1}^{k}f_{n}(A),\text{ for all }A\in \mathcal{H}(X)
\end{equation*}%
satisfies%
\begin{equation*}
H\left( T\left( A\right) ,T\left( B\right) \right) e^{H\left( T\left(
A\right) ,T\left( B\right) \right) -H\left( A,B\right) }\leq e^{-\tau
}H\left( A,B\right) ,
\end{equation*}%
for all $A,B\in \mathcal{H}(X),$ $H\left( T\left( A\right) ,T\left( B\right)
\right) \neq 0$. Using Theorem 2.1, the result follows. $\square $\medskip

\noindent \textbf{Corollary 2.6.} \ \ Let $(X,d)$ be a complete $b$-metric
space and $(X;f_{n},n=1,2,$\textperiodcentered \textperiodcentered
\textperiodcentered $,k)$ be iterated function system. Suppose that each $%
f_{i}$ for $i=1,2,...,k$ is a mapping on $X$ satisfying%
\begin{equation*}
d\left( f_{i}x,f_{i}y\right) (d\left( f_{i}x,f_{i}y\right) +1)\leq e^{-\tau
\left( d\left( x,y\right) \right) }d\left( x,y\right) (d\left( x,y\right)
+1),
\end{equation*}%
for all $x,y\in X,$ $f_{i}x\neq f_{i}y,$ where $\tau \in \Upsilon .$ Then
the mapping $T:\mathcal{H}(X)\rightarrow \mathcal{H}(X)$ defined in Theorem
2.1 has a unique fixed point in $\mathcal{H}\left( X\right) .$ Furthermore,
for any set $A_{0}\in \mathcal{H}\left( X\right) $, the sequence of compact
sets $\{A_{0},T\left( A_{0}\right) ,T^{2}\left( A_{0}\right) ,...\}$
converges to a fixed point of $T$.

\noindent \textit{Proof.} \ \ \ By taking $F\left( \lambda \right) =\ln
\left( \lambda ^{2}+\lambda \right) +\lambda ,$ $\lambda >0$ in Theorem
1.10, we obtain that each mapping $f_{i}$ for $i=1,2,...,k$ on $X\ $satisfies%
\begin{equation*}
d\left( f_{i}x,f_{i}y\right) (d\left( f_{i}x,f_{i}y\right) +1)\leq e^{-\tau
\left( d\left( x,y\right) \right) }d\left( x,y\right) (d\left( x,y\right)
+1),
\end{equation*}%
for all $x,y\in X,$ $f_{i}x\neq f_{i}y,$ where $\tau \in \Upsilon .$ Again
it follows from Theorem 1.10 that the mapping $T:\mathcal{H}(X)\rightarrow 
\mathcal{H}(X)$ defined by%
\begin{equation*}
T(A)=\cup _{n=1}^{k}f_{n}(A),\text{ for all }A\in \mathcal{H}(X)
\end{equation*}%
satisfies%
\begin{equation*}
H\left( T\left( A\right) ,T\left( B\right) \right) (H\left( T\left( A\right)
,T\left( B\right) \right) +1)\leq e^{-\tau \left( H(A,B)\right) }H\left(
A,B\right) (H\left( A,B\right) +1),
\end{equation*}%
for all $A,B\in \mathcal{H}(X),$ $H\left( T\left( A\right) ,T\left( B\right)
\right) \neq 0$. Using Theorem 2.1, the result follows. $\square $\medskip

\noindent \textbf{Corollary 2.7.} \ \ Let $(X,d)$ be a complete $b$-metric
space and $(X;f_{n},n=1,2,$\textperiodcentered \textperiodcentered
\textperiodcentered $,k)$ be iterated function system. Suppose that each $%
f_{i}$ for $i=1,2,...,k$ is a mapping on $X$ satisfying%
\begin{equation*}
d\left( f_{i}x,f_{i}y\right) \leq \frac{1}{(1+\tau \left( d(x,y)\right) 
\sqrt{d\left( x,y\right) }}d\left( x,y\right) ,
\end{equation*}%
for all $x,y\in X,$ $f_{i}x\neq f_{i}y,$ where $\tau \in \Upsilon .$ Then
the mapping $T:\mathcal{H}(X)\rightarrow \mathcal{H}(X)$ defined in Theorem
2.1 has a unique fixed point $\mathcal{H}\left( X\right) .$ Furthermore, for
any set $A_{0}\in \mathcal{H}\left( X\right) $, the sequence of compact sets 
$\{A_{0},T\left( A_{0}\right) ,T^{2}\left( A_{0}\right) ,...\}$ converges to
a fixed point of $T$.

\noindent \textit{Proof.} \ \ \ Take $F\left( \lambda \right) =-1/\sqrt{%
\lambda },$ $\lambda >0$ in Theorem 1.10, then each mapping $f_{i}$ for $%
i=1,2,...,k$ on $X\ $satisfies%
\begin{equation*}
d\left( f_{i}x,f_{i}y\right) \leq \frac{1}{(1+\tau \left( d(x,y)\right) 
\sqrt{d\left( x,y\right) })^{2}}d\left( x,y\right) ,\text{ for all }x,y\in X,%
\text{ }f_{i}x\neq f_{i}y,
\end{equation*}%
where $\tau \in \Upsilon .$ Again it follows form Theorem 1.10 that the
mapping $T:\mathcal{H}(X)\rightarrow \mathcal{H}(X)$ defined by%
\begin{equation*}
T(A)=\cup _{n=1}^{k}f_{n}(A),\text{ for all }A\in \mathcal{H}(X)
\end{equation*}%
satisfies%
\begin{equation*}
H\left( T\left( A\right) ,T\left( B\right) \right) \leq \frac{1}{(1+\tau
\left( H\left( A,B\right) \right) \sqrt{H\left( A,B\right) })^{2}}H\left(
A,B\right) ,
\end{equation*}%
for all $A,B\in \mathcal{H}(X),$ $H\left( T\left( A\right) ,T\left( B\right)
\right) \neq 0$. Using Theorem 2.1, the result follows. $\square $\medskip

\noindent \textbf{Acknowledgement}. The first author is grateful to the
Erasmus Mundus project FUSION for supporting the research visit to M\"{a}%
lardalen University, Sweden and to the Division of Applied Mathematics at
the School of Education, Culture and Communication for creating excellent
research environment.\newline

\end{document}